\documentclass[11pt]{article}

\usepackage{amssymb,amsmath,amsthm}
\usepackage{epsfig}
\usepackage{psfrag}
\newtheorem{Theorem}{\sc Theorem}

\begin{document}

\title{Picture invariants and the isomorphism problem for complex semisimple 
Lie 
algebras}
\author{Vijay Kodiyalam and K. N. Raghavan\\
Institute of Mathematical Sciences, Chennai,\\
INDIA}

\maketitle

\begin{abstract}
Consider complex semisimple Lie algebras of a given dimension specified
by their structure constants.
We describe a finite collection of rational functions in the structure
constants that form a complete set of invariants: two sets of structure
constants determine isomorphic algebras if and only if each of these
rational functions has the same value on the two sets.
This collection is indexed by all chord diagrams with the number
of chords being bounded by a function of the dimension.
\end{abstract}

Our goal in this paper is to explicitly write down a finite list of rational 
functions
in the structure constants of a complex semisimple Lie algebra of dimension $n$
such that two such Lie algebras are isomorphic if and only if the values of 
these rational functions agree for the two sets of structure constants.
Rather surprisingly - at least to us - follows the fact that the weight 
systems associated to the Lie algebra and its adjoint representation - 
see \cite{Brn} for the precise definitions - 
are a complete invariant.

We will begin by sketching the proof of our result. Semisimple
Lie algebra structures on a fixed complex vector space are seen to form  
an affine variety. This variety admits a natural action by a general linear 
group
with orbits corresponding to isomorphism classes. The basic Lie algebra
cohomology vanishing result implies that all orbits are open - by
a standard deformation theory argument - and hence that there are only
finitely many orbits which are all closed as well. Computational invariant 
theory
implies that such orbits may be distinguished by invariant polynomials in the 
coordinate ring of the variety
with explicitly bounded degree. 
Finally, an explicit pictorial spanning set of the invariant subring that is 
presented 
in \cite{DttKdySnd} is used together with graphical arguments to show
that chord diagrams specify invariants that form a complete set.
The transition from polynomial to rational occurs because the space of 
semisimple Lie algebras of dimension $n$ is only locally closed in $n^3$
dimensional affine space.
We note that the structure of this proof basically mimics that in 
\cite{DttKdySnd}.

For the rest of this paper, fix a positive integer $n$ wich will be the 
dimension of a complex vector space $V$ with a distinguished basis 
$v_1,v_2,\cdots,v_n$.
A Lie algebra structure on $V$ is defined by a set of structure constants
$\mu_{ij}^k$ where $\mu : V \otimes V \rightarrow V$ is the bracket.
Here, and in the sequel, all indices such as $i,j,k$ will range from $1$ to 
$n$ and we will consistently use the summation convention. Thus
explicitly, $\mu(v_i \otimes v_j) = \mu_{ij}^kv_k$ (where the summation over 
$k$ is suppressed).
The set of Lie algebra structures on $V$ is naturally an
affine variety in $V^* \otimes V^* \otimes V$ specified by the two tensor 
equations ensuring that the bracket is antisymmetric and satisfies the Jacobi
identity.

\begin{figure}[!h]
\begin{center}
\psfrag{M}{\huge $\mu$}
\psfrag{=}{\huge $=$}
\psfrag{0}{\huge $0$}
\psfrag{+}{\huge $+$}
\resizebox{6.0cm}{!}{\includegraphics{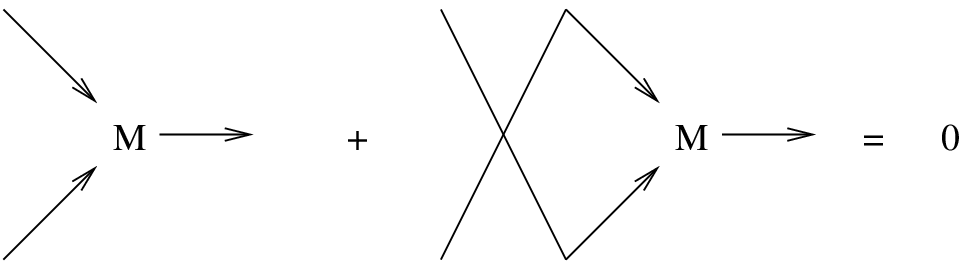}}
\end{center}
\caption{Antisymmetry of the bracket}\label{fig:antisymmetry}
\end{figure}

\begin{figure}[!h]
\begin{center}
\psfrag{M}{\huge $\mu$}
\psfrag{=}{\huge $=$}
\psfrag{0}{\huge $0$}
\psfrag{+}{\huge $+$}
\resizebox{12.0cm}{!}{\includegraphics{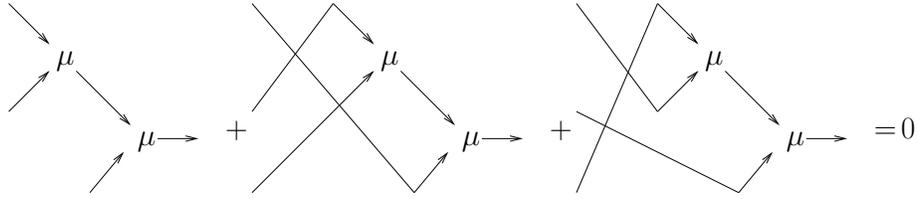}}
\end{center}
\caption{The Jacobi identity}\label{fig:jacobi}
\end{figure}

It will be important for us to use the pictorial notation for tensors due to
Kuperberg \cite{Kpr} as explained in Kauffman and Radford \cite{KffRdf}.
We will give a very brief summary of this as applied to our situation. The Lie 
algebra structure map
is represented by the picture:
$$
\begin{array}{ccc}
\searrow & & \\ & \mu & \rightarrow \\ \nearrow & &
\end{array}\ \ \ \ \ \ 
$$
while the equations defining the variety of Lie algebra structures are shown in
Figures \ref{fig:antisymmetry} and \ref{fig:jacobi}.
By convention, the inputs for a picture are read anticlockwise starting from 
the output and the 
outputs are read clockwise starting from the input. (Since this is a confusing 
point, we remark 
that this convention makes sense only when our pictures have at least one 
input and output. Otherwise we will have to explicitly indicate which is the
first input and which the first output of a picture.) Pictorial equations have 
translations in terms
of structure constants and so the reader who is suspicious of pictorial
arguments may translate everything into the language of structure constants.

\begin{figure}[!h]
\begin{center}
\psfrag{M}{\huge $\mu$}
\psfrag{B}{\huge $B$}
\psfrag{=}{\huge $=$}
\psfrag{0}{\huge $0$}
\psfrag{+}{\huge $+$}
\resizebox{5.0cm}{!}{\includegraphics{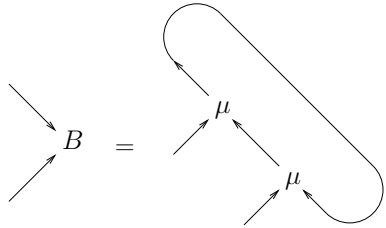}}
\end{center}
\caption{The Killing form}\label{fig:cartankilling}
\end{figure}

\begin{figure}[!h]
\begin{center}
\psfrag{M}{\huge $\mu$}
\psfrag{T}{\huge $\theta$}
\psfrag{=}{\huge $=$}
\psfrag{0}{\huge $0$}
\psfrag{+}{\huge $+$}
\resizebox{5.0cm}{!}{\includegraphics{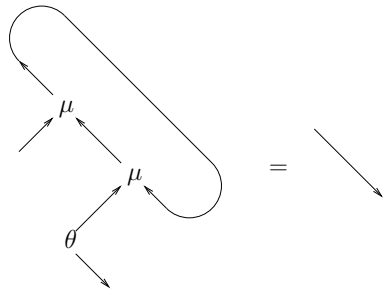}}
\end{center}
\caption{Non-degeneracy of the Killing form}\label{fig:casimir}
\end{figure}

In the variety of Lie algebra structures, those that are semisimple form an
open set defined by the non-degeneracy of the Killing form. This 
symmetric form,
denoted $B: V \otimes V \rightarrow {\mathbb C}$ is represented as in 
Figure~\ref{fig:cartankilling}.
To say that this form is non-degenerate implies the existence of an element
$\theta \in V \otimes V$ (neccessarily unique and symmetric) 
satisfying the equation in Figure \ref{fig:casimir}.
We recall here that the picture $\rightarrow$ represents the identity 
endomorphism of $V$.
Also note that since $\theta$ is symmetric, it is unneccesary to specify the 
order
of its outputs.

Since we wish to deal with varieties rather than open subsets of such,
we consider the subvariety $X$ in $(V^* \otimes V^* \otimes V) \oplus (V 
\otimes V)$
with coordinate functions $\mu_{ij}^k$ and $\theta^{ij}$ defined by the
equations in Figures \ref{fig:antisymmetry}, \ref{fig:jacobi} and 
\ref{fig:casimir}.
It should be clear that points of $X$ are in bijective correspondence with 
semisimple Lie algebra structures on $V$ and that the $\mu_{ij}^k$ components
of a point of $X$ determine - via the equation in Figure \ref{fig:casimir} - 
its $\theta^{ij}$ components.

The natural action of the group $G = GL(V)$ on $(V^* \otimes V^* \otimes V) 
\oplus 
(V \otimes V)$ restricts to an action on $X$ and the orbits correspond to
isomorphism classes of $n$-dimensional complex semsimple Lie algebras. We 
assert that all $GL(V)$ orbits are open in $X$. While this is the well known
rigidity of semisimple Lie algebras (see \cite{Grs} or \cite{NjnRch}), for 
completeness, we sketch a proof.

\begin{figure}[!h]
\begin{center}
\psfrag{M}{\huge $\mu$}
\psfrag{N}{\huge $\mu^\prime$}
\psfrag{=}{\huge $=$}
\psfrag{0}{\huge $0$}
\psfrag{+}{\huge $+$}
\resizebox{12.0cm}{!}{\includegraphics{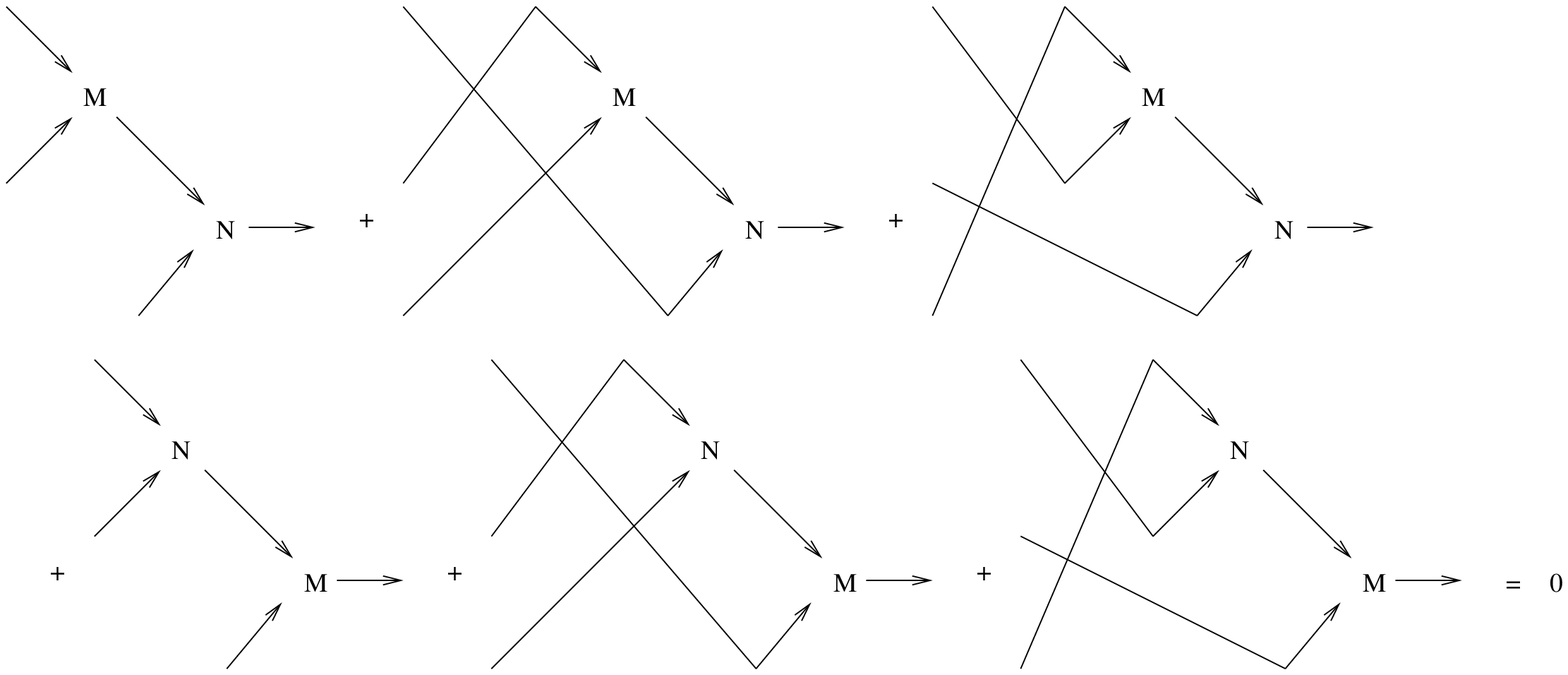}}
\end{center}
\caption{Derivative of the Jacobi identity}\label{fig:diffjacobi}
\end{figure}

\begin{figure}[!h]
\begin{center}
\psfrag{M}{\huge $\mu$}
\psfrag{N}{\huge $\mu^\prime$}
\psfrag{T}{\huge $\theta$}
\psfrag{S}{\huge $\theta^ \prime$}
\psfrag{=}{\huge $=$}
\psfrag{0}{\huge $0$}
\psfrag{+}{\huge $+$}
\resizebox{10.0cm}{!}{\includegraphics{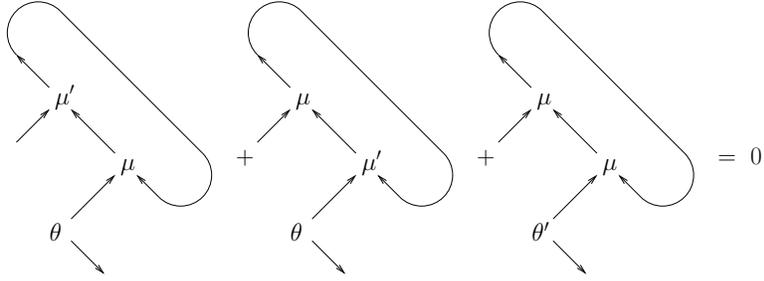}}
\end{center}
\caption{Derivative of the non-degeneracy equation}\label{fig:diffcasimir}
\end{figure}

From the equations defining the variety $X$, it follows that the
Zariski tangent space of $X$ at a point $(\mu,\theta)$ is contained in the
solution space of the linear equations in $(\mu^\prime,\theta^\prime)$ given
in Figures \ref{fig:diffantisymmetry}, \ref{fig:diffjacobi} 
and~\ref{fig:diffcasimir}.
\begin{figure}[!h]
\begin{center}
\psfrag{M}{\huge $\mu^\prime$}
\psfrag{=}{\huge $=$}
\psfrag{0}{\huge $0$}
\psfrag{+}{\huge $+$}
\resizebox{6.0cm}{!}{\includegraphics{antisymmetry.eps}}
\end{center}
\caption{Derivative of the antisymmetry equation}\label{fig:diffantisymmetry}
\end{figure}

\begin{figure}[!h]
\begin{center}
\psfrag{M}{\huge $\mu$}
\psfrag{N}{\huge $\mu^\prime$}
\psfrag{T}{\huge $\theta$}
\psfrag{S}{\huge $\theta^ \prime$}
\psfrag{=}{\huge $=$}
\psfrag{0}{\huge $0$}
\psfrag{+}{\huge $+$}
\resizebox{7.0cm}{!}{\includegraphics{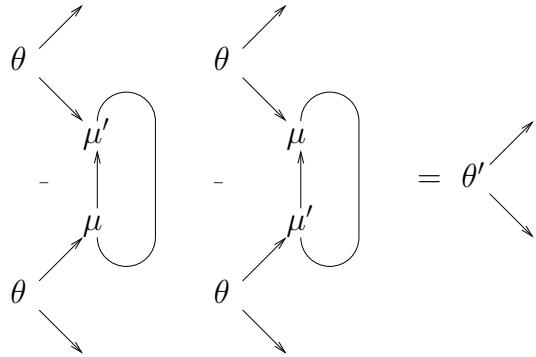}}
\end{center}
\caption{Solving for $\theta^\prime$}\label{fig:solve}
\end{figure}

We may use the equations in Figures \ref{fig:casimir} and 
\ref{fig:diffcasimir} and the symmetry of $B$ and $\theta$ to ``solve" for
$\theta^\prime$ in terms of $\mu^\prime, \mu$ and $\theta$ as in Figure 
\ref{fig:solve}.
Thus the Zariski tangent space is identified with a subspace of the solution 
space 
of the linear equations in Figures \ref{fig:diffantisymmetry} and 
\ref{fig:diffjacobi} in the variables $\mu^\prime$.

\begin{figure}[!h]
\begin{center}
\psfrag{M}{\huge $\mu$}
\psfrag{Mt}{\huge $\mu_t$}
\psfrag{T}{\huge $\theta$}
\psfrag{Tt}{\huge $\theta_t$}
\psfrag{=}{\huge $=$}
\psfrag{0}{\huge $0$}
\psfrag{eat}{\huge $e^{At}$}
\psfrag{ebt}{\huge $e^{-At}$}
\psfrag{,}{\huge $,$}
\resizebox{10.0cm}{!}{\includegraphics{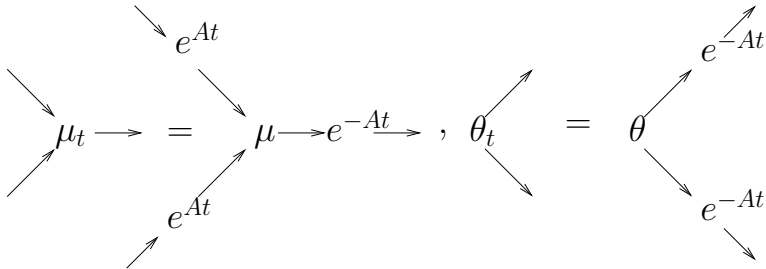}}
\end{center}
\caption{The action of $G$}\label{fig:action}
\end{figure}

\begin{figure}[!h]
\begin{center}
\psfrag{M}{\huge $\mu$}
\psfrag{N}{\huge $\mu^\prime$}
\psfrag{T}{\huge $\theta$}
\psfrag{S}{\huge $\theta^ \prime$}
\psfrag{=}{\huge $=$}
\psfrag{0}{\huge $0$}
\psfrag{+}{\huge $+$}
\psfrag{A}{\Large $A$}
\resizebox{9.0cm}{!}{\includegraphics{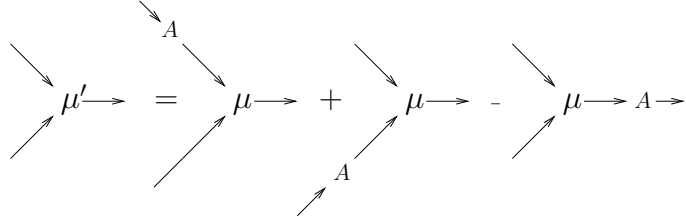}}
\end{center}
\caption{An obvious solution}\label{fig:solution}
\end{figure}

Some ``obvious elements" in the tangent space are obtained from the $G$-action
by pushing forward the Lie algebra ${\mathfrak {gl}}(V)$ of $G$ by the map 
induced by:
$G \ni T \mapsto T.(\mu,\theta) \in X$. 
Explicitly, consider $A \in {\mathfrak {gl}}(V)$. The one-parameter subgroup 
$e^{At}$
acts on $(\mu,\theta) = (\mu_0,\theta_0)$ as illustrated in Figure 
\ref{fig:action}.
Taking the derivative with respect to $t$ and evaluating at $t=0$ gives
the pictorial equation in Figure \ref{fig:solution}
as an obvious solution to the equations in Figures \ref{fig:diffantisymmetry} 
and 
\ref{fig:diffjacobi}.

Comparing with the complex below - see \S 7.2.1 of \cite{GdmWll} -
$$
0 \rightarrow {\mathfrak g} \rightarrow Hom({\mathfrak g},{\mathfrak g}) 
\rightarrow 
Hom(\wedge^2 {\mathfrak g},{\mathfrak g}) \rightarrow Hom(\wedge^3 {\mathfrak 
g},{\mathfrak 
g}) \rightarrow \cdots
$$ 
that calculates the cohomology of the adjoint module ${\mathfrak g}$ over the 
Lie 
algebra ${\mathfrak g}$ corresponding to the point $(\mu,\theta)$  of $X$,
we see that the equations in Figures \ref{fig:diffantisymmetry} and  
\ref{fig:diffjacobi}  exactly translate to saying that $\mu^\prime$ is in 
$Hom(\wedge^2 {\mathfrak g},{\mathfrak g})$ and is a cocycle in 
the complex while the obvious solutions are the coboundaries.

Since ${\mathfrak g}$ is semisimple, this cohomology vanishes - see Theorem 
3.12.1 of \cite{Vrd} -  which may be 
interpreted to mean that the map $G \rightarrow X$ considered above
is surjective at the tangent space level at the identity of $G$.
Since this holds for each point in the orbit, a simple geometric argument - 
see Corollary 1.5 of \cite{Stf} - shows that the orbit is Zariski open in $X$.

Hence there are only finitely many orbits each of which is also Zariski closed
(and the variety $X$ is non-singular).
We now appeal to what \cite{MmfFgrKrw} refer to as the `only really important 
geometric
property implied by the reductivity of $G$' - see Corollary 1.2 of Chapter 1, 
\S2  - any two orbits are separated by an element in $R^G$ where $R = {\mathbb 
C}[\mu_{ij}^k,\theta^{ij}]$ is the coordinate ring of the ambient space. It
remains to describe this invariant subring explicitly. This is exactly what
the main result of \cite{DttKdySnd} does. According to that, a linear spanning
set of the invariant ring is given by all ``closed pictures" that can be formed
using the two tensors $\mu$ and $\theta$. Further, from computational 
invariant 
theory - see  \S4.7 of \cite{DrkKmp} - $R^G$ is generated as an algebra by its 
elements of degree at most
$\frac{3}{8} (n^3 + n^2) (n + 1)^2 (2n + 1)^{2n^2}$ - the numbers
$n^3 + n^2$, $n+1$, $2n +1$ and $n^2$ being upper bounds for what they call
$r,C,A$ and $m$ respectively. 
Since any closed picture made up of $\mu$'s and $\theta$'s has twice as many 
$\mu$'s as $\theta$'s it follows that 
closed pictures involving at most $k(n) = \frac{1}{8} (n^3 + n^2) (n + 1)^2 
(2n + 1)^{2n^2}$ 
$\theta$'s separate isomorphism classes of complex semisimple
Lie algebras of dimension $n$.

We will now analyse an arbitrary closed picture involving $\mu$'s and 
$\theta$'s and show that it can be reduced modulo the relations in Figures
\ref{fig:antisymmetry}, \ref{fig:jacobi} and \ref{fig:casimir} to linear 
combinations of
products of 
chord diagrams. The interpretation of a chord diagram as a
closed picture made out of the tensors $\mu$ and $\theta$ is most easily 
illustrated
by an example as in Figure \ref{fig:chord} below. Briefly, 
\begin{figure}[!h]
\begin{center}
\psfrag{M}{\huge $\mu$}
\psfrag{N}{\huge $\mu^\prime$}
\psfrag{T}{\huge $\theta$}
\psfrag{S}{\huge $\theta^ \prime$}
\psfrag{=}{\huge $=$}
\psfrag{0}{\huge $0$}
\psfrag{+}{\huge $+$}
\resizebox{12.0cm}{!}{\includegraphics{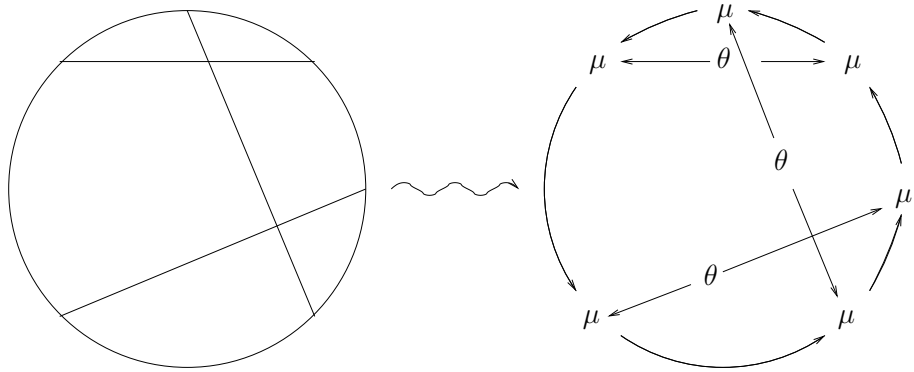}}
\end{center}
\caption{Tensor interpretation of chord diagrams}\label{fig:chord}
\end{figure}
put in a $\mu$ at each point of the circle determined by the chords, put in a 
$\theta$ on each of the chords and orient arrows as
in the example.

\begin{figure}[!h]
\begin{center}
\psfrag{M}{\huge $\mu$}
\psfrag{N}{\huge $\mu^\prime$}
\psfrag{T}{\huge $\theta$}
\psfrag{S}{\huge $\theta^ \prime$}
\psfrag{=}{\huge $=$}
\psfrag{0}{\huge $0$}
\psfrag{+}{\huge $+$}
\resizebox{6.0cm}{!}{\includegraphics{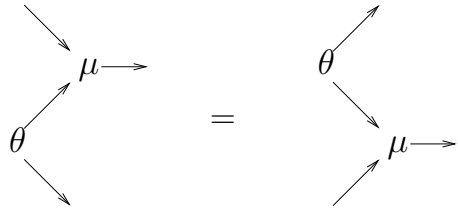}}
\end{center}
\caption{A useful pictorial equality}\label{fig:picequality}
\end{figure}

\begin{figure}[!h]
\begin{center}
\psfrag{M}{\huge $\mu$}
\psfrag{B}{\Large $B$}
\psfrag{N}{\huge $\mu^\prime$}
\psfrag{T}{\huge $\theta$}
\psfrag{S}{\huge $\theta^ \prime$}
\psfrag{=}{\huge $=$}
\psfrag{0}{\huge $0$}
\psfrag{+}{\huge $+$}
\psfrag{*}{\huge $*$}
\resizebox{6.0cm}{!}{\includegraphics{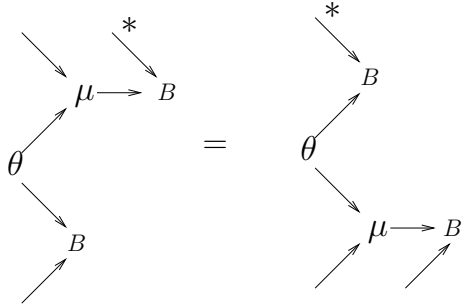}}
\end{center}
\caption{Proving the equality of Figure \ref{fig:picequality}}\label{fig:ets}
\end{figure}

\begin{figure}[!h]
\begin{center}
\psfrag{M}{\huge $\mu$}
\psfrag{N}{\huge $\mu^\prime$}
\psfrag{T}{\huge $\theta$}
\psfrag{S}{\huge $\theta^ \prime$}
\psfrag{=}{\huge $=$}
\psfrag{0}{\huge $0$}
\psfrag{+}{\huge $+$}
\resizebox{7.0cm}{!}{\includegraphics{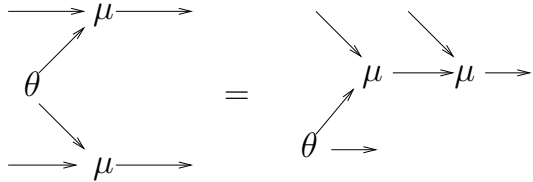}}
\end{center}
\caption{A corollary}\label{fig:corollary}
\end{figure}

Before proceeding with this proof, we pause to observe a useful
pictorial consequence of the cyclic symmetry of the trilinear form:
$a \otimes b \otimes c \mapsto B([a,b],c)$. This is the pictorial
equality in Figure \ref{fig:picequality}.
To prove this equation, use the non-degeneracy of $B$ to reduce to proving the
equality in Figure \ref{fig:ets}. (The $*$'s are supposed to indicate which
is the first input since these are pictures without outputs.)
But now, using the symmetry of $B$ and $\theta$ and the equation in Figure 
\ref{fig:casimir}, we see that this is pictorial translation of the cyclic
symmetry of the trilinear form above.
What we will actually use is the corollary in Figure \ref{fig:corollary} of the
equation in Figure \ref{fig:picequality}, the easy proof of which is
left to the reader.

Consider now an arbitrary closed picture, say $P$, involving $\mu$'s and 
$\theta$'s. 
We will first discuss the case that $P$ is connected.
Suppose that the number of $\theta$'s in $P$ is $k$ so that the number of 
$\mu$'s is $2k$.
Removing all $\theta$'s from the picture gives a picture, say $Q$, that 
specifies a map 
from $V^{\otimes 2k}$ to
${\mathbb C}$ (after choosing an order on its inputs) that is expressible 
using only $\mu$'s. While $Q$ is not 
neccessarily connected, note that in the equation in Figure
\ref{fig:corollary}, the picture on the right remains connected after removing
$\theta$ while that on the left is disconnected. Thus by using this corollary
as often as needed we may asume that $Q$ is connected.

Consider $Q$ as a connected directed graph where each vertex
is either univalent with outdegree $1$ or trivalent with indegree $2$ and
outdegree $1$.
We claim that there is exactly one directed cycle in this graph.
To see this, observe that $Q$ has $4k$ vertices - $2k$ of them trivalent and
$2k$ univalent - and $4k$ edges - by counting valencies. So, regarded
as an undirected graph, it has exactly one cycle.
To see that this cycle is, in fact, directed, consider it as a subgraph of
the original directed graph and look at the list of pairs $(i,o)$ where
$i$ is the indegree and $o$ is the outdegree - in the cycle - of each
vertex. 
Each pair is one of $(0,2),(2,0),(1,1)$ and the sum of the $i$'s is equal
to the sum of the $o$'s.
Hence the number of vertices in this cycle with data $(0,2)$ is equal
to that with data $(2,0)$.
Since each vertex is a $\mu$, it is not possible that $(0,2)$ occurs for
any vertex and therefore each vertex is of type $(1,1)$.
This means that the cycle is directed.

\begin{figure}[!h]
\begin{center}
\psfrag{M}{\huge $\mu$}
\psfrag{N}{\huge $\mu^\prime$}
\psfrag{T}{\huge $\theta$}
\psfrag{S}{\huge $\theta^ \prime$}
\psfrag{p1}{\huge $P_1$}
\psfrag{p2}{\huge $P_2$}
\psfrag{p3}{\huge $P_3$}
\psfrag{pi}{\huge $P_i$}
\psfrag{ptm1}{\huge $P_{t-1}$}
\psfrag{pt}{\huge $P_t$}
\psfrag{0}{\huge $0$}
\psfrag{+}{\huge $+$}
\resizebox{5.0cm}{!}{\includegraphics{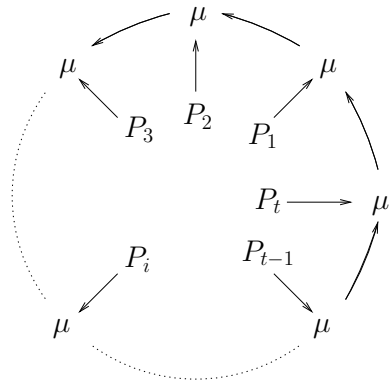}}
\end{center}
\caption{Reduction}\label{fig:reduction}
\end{figure}

Next, note that by using the antisymmetry relation we may assume that for all 
the $\mu$'s occuring in the cycle, the
in-edge of $\mu$ that is part of the cycle is the one
corresponding to its second input. Thus we assume that $Q$ is of the form 
shown in Figure \ref{fig:reduction} where the $P_i$ 
are pictures
representing maps from some tensor power of $V$ to $V$ that are made up of 
$\mu$'s and do not have cycles.

\begin{figure}[!h]
\begin{center}
\psfrag{M}{\huge $\mu$}
\psfrag{=}{\huge $=$}
\psfrag{-}{\huge $-$}
\psfrag{N}{\huge $\mu^\prime$}
\psfrag{T}{\huge $\theta$}
\psfrag{S}{\huge $\theta^ \prime$}
\psfrag{p1}{\huge $P_1$}
\psfrag{p2}{\huge $P_2$}
\psfrag{p3}{\huge $P_3$}
\psfrag{pi}{\huge $P_i$}
\psfrag{ptm1}{\huge $P_{t-1}$}
\psfrag{pt}{\huge $P_t$}
\psfrag{0}{\huge $0$}
\psfrag{+}{\huge $+$}
\resizebox{10.0cm}{!}{\includegraphics{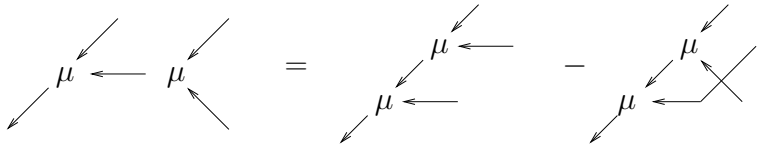}}
\end{center}
\caption{The Jacobi identity}\label{fig:anotherjacobi}
\end{figure}

Next we will induce on the total number of $\mu$'s occuring in all the $P_i$
to show that we may reduce to the case that this number is $0$. This is easy
and follows immediately from the Jacobi identity written in the form in Figure
\ref{fig:anotherjacobi}. 
Basically, we may get rid of the complexity inside
a circle at the expense of increasing its size.

We have therefore seen that modulo the relations in Figures 
\ref{fig:antisymmetry}, \ref{fig:jacobi} and \ref{fig:casimir} any connected 
closed
picture made up of the tensors $\mu$ and $\theta$ is equivalent to a
linear combination of 
%products of 
pictures associated to chord diagrams.
In case the original picture is disconnected, we will have to allow for
a linear combination of 
products of 
pictures associated to chord diagrams.
A review of the proof shows that if the original picture
involved $k$ $\theta$'s, then so do the resulting chord diagrams. 

While a closed picture made up of $\theta$'s and $\mu$'s corresponds to an 
element in the polynomial ring
$R = {\mathbb 
C}[\mu_{ij}^k,\theta^{ij}]$, recall that for a point lying on $X$, 
$\theta^{ij}$
may be expressed as a rational function of the $\mu_{ij}^k$ and thus each
such closed picture - in particular, a chord diagram - determines a rational 
function in the $\mu_{ij}^k$ on $X$.

To summarise, we have proved the following theorem.

\begin{Theorem}
Any chord diagram determines a rational function of $n^3$ variables
$\mu^k_{ij}$ for any positive integer $n$.
The values of these rational functions at the structure constants of
a complex semisimple Lie algebra of dimension $n$ are isomorphism invariants.
If, for two complex semisimple Lie algebras of dimension $n$, these values 
agree for all chord diagrams with at most $k(n) = \frac{1}{8} (n^3 + n^2) 
(n + 1)^2 (2n + 1)^{2n^2}$ chords, then the two
Lie algebras are isomorphic.\qed
\end{Theorem}

We conclude by observing that these numbers associated to a chord diagram
and a complex semisimple Lie algebra are exactly those studied in \cite{Brn}
and are essentially equivalent to the weight systems associated to
such a Lie algebra and its adjoint representation.

\medskip

\noindent{\bf Acknowledgements.} VK would like to thank Murray Gerstenhaber and
Willem de Graaf for
 some helpful references and e-mails.

\end{document}